\documentclass[10pt]{article}
\usepackage{latexsym}
\usepackage{amssymb}
\newtheorem{thm}{Theorem}

\newcommand{\Z}{\hbox{\bf Z}}

\newcommand{\beeq}{\begin{eqnarray*}}
\newcommand{\eneq}{\end{eqnarray*}}
\newcommand{\proof}{\noindent {\it Proof.\hspace{4mm}}}

\newcommand{\qfd}{\hfill $\fbox{}$\vspace{4mm}}
\def\newpic#1{%
\def\emline##1##2##3##4##5##6{%
\put(##1,##2){\special{em:point #1##3}}%
\put(##4,##5){\special{em:point #1##6}}%
\special{em:line #1##3,#1##6}}}
\newpic{}
\def\emline#1#2#3#4#5#6{%
\put(#1,#2){\special{em:moveto}}%
\put(#4,#5){\special{em:lineto}}}
\def\newpic#1{}
\title{Pappus-Desargues digraph confrontation}%Pappus vs. Desargues reattachment graphs}
\author{Italo J. Dejter
\\ University of Puerto Rico \\ Rio Piedras, PR 00936-8377 \\ italo.dejter@gmail.com}
\date{}
\begin{document}
\maketitle
\begin{abstract}
\noindent Like the Coxeter graph became reattached into the Klein
graph in \cite{DCox}, the Levi graphs of the $9_3$ and $10_3$
self-dual configurations, known as the Pappus and Desargues
($k$-transitive) graphs $\mathcal P$ and $\mathcal D$ (where $k=3$),
also admit reattachments of the distance-$(k-1)$ graphs of half of
their oriented shortest cycles via orientation assignments on their
common $(k-1)$-arcs, concurrent for ${\mathcal P}$ and opposite for
$\mathcal D$, now into 2 disjoint copies of their corresponding
Menger graphs. Here, $\mathcal P$ is the unique cubic
distance-transitive (or CDT) graph with the concurrent-reattachment
behavior while $\mathcal D$ is one of 7 CDT graphs with the
opposite-reattachment behavior, that include the Coxeter graph.
Thus, $\mathcal P$ and $\mathcal D$ confront each other in these
respects, obtained via $\mathcal C$-ultrahomogeneous graph
techniques \cite{orestes,I} that allow to characterize the obtained
reattachment Menger graphs in the same terms.
\end{abstract}

%\noindent{\bf Keywords:} ultrahomogeneous graph; digraph; shortest cycle

%\noindent{\bf 2000 Mathematics subject classification:} 05C62, 05B30, 05C20, 05C38

\section{Preliminaries}

\noindent Given a collection $\mathcal C$ of (di)graphs closed under
isomorphisms, a (di)graph $G$ is said to be $\mathcal C$-{\it
ul\-tra\-ho\-mo\-ge\-neous} (or $\mathcal C$-{\it UH})
\cite{orestes,I} if every isomorphism between 2 induced members of
$\mathcal C$ in $G$ extends to an auto\-mor\-phism of $G$. If
${\mathcal C}=\{H\}$ is the isomorphism class of a (di)graph $H$, we
say that such a $G$ is $\{H\}$-{\it UH} (or $H$-{\it UH}). In
\cite{I}, $\mathcal C$-UH graphs are studied when ${\mathcal C}$ is
the collection of either {\bf(a)} the complete graphs, or {\bf(b)}
the disjoint unions of complete graphs, or {\bf(c)} the complements
of those unions.\bigskip

\noindent We consider any undirected graph $G$ as a digraph by
taking each edge $e$ of $G$ as a pair of oppositely oriented (or
{\it O-O}) arcs $\vec{e}$ and $(\vec{e})^{-1}$. Then {\it cohering})
(or {\it fastening}, or {\it zipping}) $\vec{e}$ and
$(\vec{e})^{-1}$ (meaning that we take the union of $\vec{e}$ and
$(\vec{e})^{-1}$) allows to obtain precisely $e$, a simple technique
to be used below. In other words, $G$ is a graph taken as a digraph,
that is, for any $2$ adjacent vertices $u,v\in V(G)$, the arcs
$\vec{e}=(u,v)$ and $(\vec{e})^{-1}=(v,u)$ are both present in the
set $A(G)$ of arcs of $G$, with the union
$\vec{e}\,\cup(\vec{e})^{-1}$ interpreted as the edge $e\in E(G)$ of
$G$. If we write $\vec{f}=(\vec{e})^{-1}$, then clearly
$(\vec{f}\,)^{-1}=\vec{e}$ and $f=e$.

\subsection{Coherent $\mathcal C$-(ultra)homogeneous graphs}

\noindent Let $M$ be an induced subgraph of a graph $H$ and let $G$
be both an $M$-UH and an $H$-UH graph. We say that $G$ is an
$\{H\}_{M}$-{\it UH} graph if, for each copy $H_0$ of $H$ induced in
$G$ and containing a copy $M_0$ of $M$, there exists exactly one
copy $H_1\neq H_0$ of $H$ induced in $G$ such that $V(H_0)\cap
V(H_1)=V(M_0)$ and $E(H_0)\cap E(H_1)=E(M_0)$. These vertex and edge
conditions can be condensed as $H_0\cap H_1=M_0$. We say that such a
$G$ is {\it coherent}. This is generalized by saying that an
$\{H\}_M$-UH graph $G$ is an $\ell${\it -coherent} $\{H\}_M$-{\it
UH} graph if, given a copy $H_0$ of $H$ induced in $G$ and
containing a copy $M_0$ of $M$, there exist exactly $\ell$ copies
$H_i\neq H_0$ of $H$ induced in $G$ such that $H_i\cap H_0\supseteq
M_0$, for each $i=1,2,\ldots,\ell$, with $H_1\cap H_0=M_0$.\bigskip

\noindent If $G$
is coherent $\{H\}_M$-UH and $K$ is both subgraph of $H$ and
supergraph of $M$, we say that $G$ is $\{K\not\subset H\}_M$-UH if
every isomorphism between 2 induced copies of $K$ in $G$ not
contained in any copy of $H$ in $G$ extends to an auto\-mor\-phism
of $G$. If, under these conditions, each copy of $M$ induced in $G$
coincides with the intersection of exactly one copy of $H$ and
exactly one copy of $K\not\subset H$, then we say that $G$ is {\it
coherent} $\{H,K\}_M$-{\it UH}. This concept is used in Theorem 3
below for the Desargues graph $G=\mathcal D$.\bigskip

\noindent Let $G$ be an $M$-UH graph but not $H$-UH and assume the
isomorphism class $\mathcal H$ of $H$ in $G$ decomposes as
${\mathcal H}={\mathcal H}_0\cup{\mathcal H}_1$ so that every
isomorphism between $2$ members of ${\mathcal H}_i$ induced in $G$
extends to an auto\-mor\-phism of $G$, ($i=0,1$). If $H_i$ is a
representative of ${\mathcal H}_i$, for $i=0,1$, then we say that
$G$ is an $\{H_0,H_1\}^M$-{\it homogeneous graph} if, for each copy
$H_i$ induced in $G$ and containing a copy $M_0$ of $M$, there
exists exactly one copy $H_j$ induced in $G$, ($i,j\in\{0,1\}, i\ne
j$), with $H_i\cap H_j=M_0$. This concept is likewise extended to a
decomposition ${\mathcal H}={\mathcal H}_1\cup{\mathcal
H}_2\cup{\mathcal H}_3$ in Theorem 5, where $G={\mathcal P}$ is the
Pappus graph.

\subsection{Coherent $\mathcal C$-(ultra)homogeneous digraphs}

\noindent Let $\vec{M}$ be an induced subdigraph of a digraph
$\vec{H}$ and let $G$ be both an $\vec{M}$-UH and an $\vec{H}$-UH
digraph. We say that $G$ is an $\{\vec{H}\}_{\vec{M}}$-{\it UH}
digraph (resp. an $\{\vec{H}\}^{\vec{M}}$-{\it UH} digraph), if for
each copy $\vec{H}_0$ of $\vec{H}$ induced in $\vec{G}$ and
containing a copy $\vec{M}_0$ of $\vec{M}$, there exists exactly one
copy $\vec{H}_1\neq\vec{H}_0$ of $\vec{H}$ induced in $G$ such that
$V(\vec{H}_0)\cap V(\vec{H}_1)=V(\vec{M}_0)$ and
$A(\vec{H}_0)\cap\bar{A}(\vec{H}_1)=A(\vec{M}_0)$ (resp.
$A(\vec{H}_0)\cap{A}(\vec{H}_1)=A(\vec{M}_0)$), where
$\bar{A}(\vec{H}_1)$ is formed by those arcs $(\vec{e})^{-1}$ whose
orientations are reversed with respect to the orientations of the
arcs $\vec{e}$ of $A(\vec{H}_1)$. In either case, we may say that
such a $G$ is {\it coherent}.\bigskip

\noindent Let $G$ be an $\vec{M}$-UH graph but not $\vec{H}$-UH and
assume  that the isomorphism class $\vec{\mathcal H}$ of $\vec{H}$
in $G$ decomposes as $\vec{\mathcal H}=\vec{\mathcal
H}_0\cup\vec{\mathcal H}_1$ so that every isomorphisms between $2$
members of $\vec{\mathcal H}_i$ induced in $G$ extends to an
auto\-mor\-phism of $G$, ($i=0,1$). If $\vec{H}_i$ is a
representative of $\vec{\mathcal H}_i$, for $i=0,1$, then we say
that $G$ is an $\{\vec{H}_0,\vec{H}_1\}^{\vec{M}}$-{\it homogeneous}
graph if, for each copy $\vec{H}_i$ induced in $G$ and containing a
copy $\vec{M}_0$ of $M$, there exists exactly one induced copy
$\vec{H_j}$ in $G$, ($i,j\in\{0,1\}, i\ne j$), with
$\vec{H}_i\cap\vec{H}_j=\vec{M}_0$.

\subsection{Strongly coherent $\mathcal C$-ultrahomogeneous graphs}

\noindent Given a finite graph $H$ and a subgraph $M$ of $H$ with
$|V(H)|>3$, we say that a graph $G$ is {\it strongly coherent} (or
{\it SC}) $\{H\}_M$-{\it UH} if there is a descending sequence of
connected subgraphs $M=M_1,M_2\ldots,M_t\equiv K_2$ such that:
{\bf(a)} $M_{i+1}$ is obtained from $M_i$ by the deletion of a
vertex, for $i=1,\ldots,t-1$ and {\bf(b)} $G$ is a
$(2^i-1)$-coherent $\{H\}_{M_i}$-UH graph, for
$i=1,\ldots,t$.\bigskip

\noindent Some parameters of $\mathcal P$ and $\mathcal D$ (see for
example \cite{BCN}) can be displayed as follows:
$$\begin{array}{l|l|l|l|l|l|l|l|l}
_{G} & _{n}  & _{d} &  _{g}  & _{k} &  _{\eta} &  _{a} & _{b} &
_{h}\\ \hline
%&&&&&&&&\\ \hline
^{\mathcal P}_{\mathcal D}&^{18}_{20}&^4_5&^6_6&^3_3&^{18}_{20}&^{216}_{240}&^1_1&^1_1\\
\end{array}$$
where $n,d,g,k,\eta$ and $a$ are respectively: order, diameter,
girth, largest $\ell$ such that $G$ is $\ell$-arc transitive, number
of $g$-cycles and number of automorphisms, with $b$ (resp. $h$) $=1$
if $G$ is bipartite (resp. hamiltonian) and $=0$ otherwise. Theorem
1 below asserts that both the Pappus graph $\mathcal P$ and the
Desargues graph $\mathcal D$ are SC $\{C_6\}_{P_3}$-UH graphs,
(which is also the case of the other 10 CDT graphs, see
\cite{DCox,orestes}).

\subsection{Plan of the subsequent sections}

\noindent Given a (di)graph $\Gamma$, the {\it distance}-$(k-1)$
({\it di}){\it graph} $\Gamma^{k-1}$ of $\Gamma$ has
$V(\Gamma^{k-1})=V(\Gamma)$ and an arc $(u,v)$ for each shortest
$(k-1)$-arc in $\Gamma$ from $u$ to $v\ne u$. If $\Gamma$ is a
cycle, then $\Gamma^2$ is said to be a {\it square}. Theorem 2 below
establishes that $\mathcal D$ is a $\{\vec{C_g}\}_{\vec{P}_k}$-UH
digraph and that $\mathcal P$ is a $\{\vec{C_g}\}^{\vec{P}_k}$-UH
digraph; it deals with just a pair of the 12 cubic
distance-transitive (or CDT) graphs treated in Theorem 3 of
\cite{orestes} and is given, together with its proof, in part for
the needs of the constructions in Sections 3-6. However, we stress
here that $\mathcal P$ is the only CDT graph that is a
$\{\vec{C}_g\}^{\vec{P}_k}$-UH digraph, while $\mathcal D$ is the
second most interesting of 7 CDT graphs $G$ that are
$\{\vec{C}_g\}_{\vec{P}_k}$-UH digraphs \cite{orestes} after the
Coxeter graph, where $g$ is the girth of $k$-transitive $G$.
(Petersen, Heawood, Foster and Biggs-Smith graphs excluded here.
$K_4$, $K_{3,3}$, the 3-cube and the dodecahedral graphs and Tutte
8-cage have either $g=2(k-1)$ or $k=2$, so the equivalent of the
composed operation (2) in Section 3 below or in Section 3 of
\cite{DCox} for the Coxeter graph is less interesting).\bigskip

\noindent In Sections 3-4, the squares of oriented cycles of
$\mathcal D$ yield a coherent

\noindent $\{K_4,K_3\}_{K_2}$-UH graph by means of the O-O 2-arcs
shared (as 2-paths) by the 6-cycle square pairs. In Theorem 3, this
is shown to be the disjoint union of $2$ copies of $L(K_5)$, the
Menger graph of the self-dual $(10_3)$-configuration \cite{Cox},
(whose Levi graph is $\mathcal D$, \cite{Cox}. Compare with
\cite{DCox}, yielding the Klein graph from the Coxeter graph. Recall
that the Menger graph $\mathcal M$ of a self-dual configuration
$\mathcal S$ has as vertices its points, with any $2$ determining an
edge of $\mathcal M$ if and only if their representative points in
$\mathcal S$ are colinear \cite{Cox}. We note that $2$ different
configurations may have the same Menger graph, unless each line of
$\mathcal S$ determines a maximal clique in $\mathcal M$, which is
the case of the coherent $\mathcal C$-UH graphs in this paper.) We
finish Section 4 noting that Theorem 3 yields an infinite nested
sequence of geometric realizations of $L(K_5)$ (or of its
complement, the Petersen graph) via taking barycenters of
participating tetrahedra as vertices of subsequent tetrahedra.
Generalizing, Theorem 4 in Section 5 asserts that for $n\geq 4$ the
line graph $L(K_n)$ is a coherent $\{K_{n-1},K_3\}_{K_2}$-UH graph
containing $n$ copies of $K_{n-1}$ and ${n\choose 3}$ copies of
$K_3$. An adaptation of the previous considerations to $\mathcal P$
makes it yield, in Theorem 5 of Section 6, $\mathcal P$ in $2$
complementary ways as the Menger graph of the self-dual
$(9_3)$-configuration (whose Levi graph is $\mathcal P$, \cite{Cox})
and as the object of application of the concepts of $\mathcal
C$-homogeneous graphs and digraphs given above.

\section{$(C_6,P_3)$-UH properties of $\mathcal P$ and $\mathcal D$}

%pdsolo, cokebus, cokebuso, powerade, ogocdtg, papudesa.

\begin{thm}
Let $G$ be either $\mathcal P$ or $\mathcal D$. Then $G$ is an SC
$\{C_6\}_{P_3}$-UH graph.
\end{thm}

\proof We must see that each of $G={\mathcal P}$ and $G={\mathcal
D}$ is a $(2^{i+1}-1)$-coherent $\{C_6\}_{P_{3-i}}$-UH graph, for
$i=0,1$. Taking into account details in the proof of Theorem 2
below, each $(2-i)$-path $P=P_{3-i}$ of $G$ is seen to be shared by
exactly $2^{i+1}$ 6-cycles of $G$, for $i=0,1$. It follows that $G$
is an SC $\{C_6\}_{P_3}$-UH graph.\qfd

\noindent In both ${\mathcal P}$ and ${\mathcal D}$, there are just
$2$ $6$-cycles shared by each $2$-path. If $G$ is a
$\{\vec{C}_6\}_{\vec{P}_3}$-UH digraph, then there is an assignment of
an orientation to each 6-cycle of $G$ so that the $2$ 6-cycles
shared by each $2$-path receive opposite orientations. We say that
such an assignment is a $\{\vec{C}_6\}_{\vec{P}_3}$-{\it O-O
assignment} (or $\{\vec{C}_6\}_{\vec{P}_3}$-{\it OOA}). The
collection of $\eta$ oriented 6-cycles corresponding to the $\eta$
6-cycles of $G$, for a particular $\{\vec{C}_6\}_{\vec{P}_3}$-OOA,
is called an $\{\eta\vec{C}_6\}_{\vec{P}_3}$-{\it OOC}. Each such
cycle is written with their successive vertices between parentheses
but without separating commas, where as usual the vertex that
succeeds the last vertex of the cycle is the first vertex. Arcs are
written $(u,v)$ and 2-arcs $(u,v,w)$. Figure 1 contains
representations of $\mathcal P$ and $\mathcal D$ using the vertex
notation in the proof of Theorem 2 (plus extra-features for
$\mathcal D$ related to Figure 2 and the treatment of $\mathcal D$
in Section 4 below).

\begin{thm}
$\mathcal D$ is a $\{\vec{C}_6\}_{\vec{P}_3}$-UH digraph but
$\mathcal P$ is a $\{\vec{C}_6\}^{\vec{P}_3}$-UH digraph.
\end{thm}

\begin{figure}[htp]
\unitlength=0.4mm \special{em:linewidth 0.4pt} \linethickness{0.4pt}
\begin{picture}(271.00,105.00)
\put(82.00,5.00){\circle{2.00}}
\put(130.00,45.00){\circle{2.00}}
\put(130.00,61.00){\circle{2.00}}
\put(82.00,101.00){\circle{2.00}}
\put(34.00,61.00){\circle{2.00}}
\put(34.00,45.00){\circle{2.00}}
\put(115.00,15.00){\circle{2.00}}
\put(49.00,15.00){\circle{2.00}}
\put(49.00,91.00){\circle{2.00}}
\put(65.00,8.00){\circle{2.00}}
\put(99.00,8.00){\circle{2.00}}
\put(37.00,77.00){\circle{2.00}}
\put(37.00,29.00){\circle{2.00}}
\put(126.00,77.00){\circle{2.00}}
\put(126.00,29.00){\circle{2.00}}
\put(65.00,98.00){\circle{2.00}}
\put(115.00,91.00){\circle{2.00}}
\put(99.00,98.00){\circle{2.00}}
\emline{66.00}{97.00}{1}{129.00}{62.00}{2}
\emline{129.00}{44.00}{3}{66.00}{9.00}{4}
\emline{49.00}{90.00}{5}{49.00}{16.00}{6}
\emline{100.00}{98.00}{7}{114.00}{92.00}{8}
\emline{130.00}{60.00}{9}{130.00}{46.00}{10}
\emline{130.00}{44.00}{11}{127.00}{30.00}{12}
\emline{126.00}{28.00}{13}{116.00}{16.00}{14}
\emline{114.00}{14.00}{15}{100.00}{8.00}{16}
\emline{98.00}{7.00}{17}{83.00}{5.00}{18}
\emline{83.00}{101.00}{19}{98.00}{99.00}{20}
\emline{81.00}{5.00}{21}{66.00}{7.00}{22}
\emline{64.00}{8.00}{23}{50.00}{14.00}{24}
\emline{48.00}{16.00}{25}{37.00}{28.00}{26}
\emline{36.00}{30.00}{27}{34.00}{44.00}{28}
\emline{34.00}{46.00}{29}{34.00}{60.00}{30}
\emline{50.00}{92.00}{31}{64.00}{98.00}{32}
\emline{66.00}{99.00}{33}{81.00}{101.00}{34}
\emline{130.00}{62.00}{35}{127.00}{76.00}{36}
\emline{126.00}{78.00}{37}{116.00}{90.00}{38}
\emline{116.00}{90.00}{39}{116.00}{90.00}{40}
\emline{34.00}{62.00}{41}{36.00}{76.00}{42}
\emline{37.00}{78.00}{43}{48.00}{90.00}{44}
\emline{98.00}{97.00}{45}{38.00}{29.00}{46}
\emline{114.00}{16.00}{47}{82.00}{100.00}{48}
\emline{35.00}{45.00}{49}{125.00}{29.00}{50}
\emline{114.00}{90.00}{51}{83.00}{6.00}{52}
\emline{98.00}{9.00}{53}{38.00}{76.00}{54}
\emline{35.00}{61.00}{55}{125.00}{77.00}{56}
\put(63.00,102.00){\makebox(0,0)[cc]{$_1$}}
\put(82.00,105.00){\makebox(0,0)[cc]{$_2$}}
\put(102.00,102.00){\makebox(0,0)[cc]{$_3$}}
\put(118.00,95.00){\makebox(0,0)[cc]{$_4$}}
\put(130.00,79.00){\makebox(0,0)[cc]{$_5$}}
\put(135.00,61.00){\makebox(0,0)[cc]{$_6$}}
\put(135.00,45.00){\makebox(0,0)[cc]{$_7$}}
\put(130.00,28.00){\makebox(0,0)[cc]{$_8$}}
\put(118.00,11.00){\makebox(0,0)[cc]{$_9$}}
\put(102.00,4.00){\makebox(0,0)[cc]{$_a$}}
\put(82.00,1.00){\makebox(0,0)[cc]{$_b$}}
\put(64.00,4.00){\makebox(0,0)[cc]{$_c$}}
\put(47.00,11.00){\makebox(0,0)[cc]{$_d$}}
\put(33.00,28.00){\makebox(0,0)[cc]{$_e$}}
\put(31.00,45.00){\makebox(0,0)[cc]{$_f$}}
\put(31.00,61.00){\makebox(0,0)[cc]{$_g$}}
\put(32.00,79.00){\makebox(0,0)[cc]{$_h$}}
\put(45.00,94.00){\makebox(0,0)[cc]{$_0$}}
\put(210.00,5.00){\circle{2.00}}
\put(226.00,5.00){\circle{2.00}}
\put(266.00,45.00){\circle{2.00}}
\put(266.00,61.00){\circle{2.00}}
\put(226.00,101.00){\circle{2.00}}
\put(210.00,101.00){\circle{2.00}}
\put(170.00,61.00){\circle{2.00}}
\put(170.00,45.00){\circle{2.00}}
\put(254.00,17.00){\circle{2.00}}
\put(182.00,17.00){\circle{2.00}}
\put(182.00,89.00){\circle{2.00}}
\put(194.00,10.00){\circle{2.00}}
\put(242.00,10.00){\circle{2.00}}
\put(173.00,77.00){\circle{2.00}}
\put(173.00,29.00){\circle{2.00}}
\put(262.00,77.00){\circle{2.00}}
\put(262.00,29.00){\circle{2.00}}
\put(194.00,97.00){\circle{2.00}}
\put(242.00,97.00){\circle{2.00}}
\put(254.00,89.00){\circle{2.00}}
\emline{227.00}{5.00}{57}{241.00}{9.00}{58}
\emline{209.00}{5.00}{59}{195.00}{9.00}{60}
\emline{243.00}{11.00}{61}{253.00}{16.00}{62}
\emline{193.00}{11.00}{63}{183.00}{16.00}{64}
\emline{255.00}{18.00}{65}{261.00}{28.00}{66}
\emline{181.00}{18.00}{67}{174.00}{28.00}{68}
\emline{172.00}{30.00}{69}{170.00}{44.00}{70}
\emline{170.00}{62.00}{71}{172.00}{76.00}{72}
\emline{174.00}{78.00}{73}{181.00}{88.00}{74}
\emline{183.00}{90.00}{75}{193.00}{96.00}{76}
\emline{195.00}{98.00}{77}{209.00}{101.00}{78}
\emline{227.00}{101.00}{79}{241.00}{98.00}{80}
\emline{243.00}{96.00}{81}{253.00}{90.00}{82}
\emline{253.00}{90.00}{83}{253.00}{90.00}{84}
\emline{255.00}{88.00}{85}{261.00}{78.00}{86}
\emline{263.00}{76.00}{87}{266.00}{62.00}{88}
\emline{266.00}{44.00}{89}{263.00}{30.00}{90}
\emline{183.00}{89.00}{91}{253.00}{89.00}{92}
\emline{211.00}{100.00}{93}{241.00}{10.00}{94}
\emline{225.00}{100.00}{95}{195.00}{11.00}{96}
\emline{171.00}{61.00}{97}{265.00}{61.00}{98}
\emline{174.00}{76.00}{99}{253.00}{18.00}{100}
\emline{183.00}{18.00}{101}{261.00}{76.00}{102}
\emline{211.00}{101.00}{103}{225.00}{101.00}{104}
\emline{266.00}{60.00}{105}{266.00}{46.00}{106}
\emline{170.00}{60.00}{107}{170.00}{46.00}{108}
\emline{211.00}{5.00}{109}{225.00}{5.00}{110}
\put(190.00,101.00){\makebox(0,0)[cc]{$_{0_0}$}}
\put(210.00,105.00){\makebox(0,0)[cc]{$_{0_1}$}}
\put(227.00,105.00){\makebox(0,0)[cc]{$_{0_2}$}}
\put(246.00,100.00){\makebox(0,0)[cc]{$_{\underline{0_3}}$}}
\put(266.00,80.00){\makebox(0,0)[cc]{$_{\underline{1_1}}$}}
\put(271.00,62.00){\makebox(0,0)[cc]{$_{1_2}$}}
\put(270.00,48.00){\makebox(0,0)[cc]{$_{\underline{1_3}}$}}
\put(268.00,28.00){\makebox(0,0)[cc]{$_{2_0}$}}
\put(260.00,16.00){\makebox(0,0)[cc]{$_{2_1}$}}
\put(209.00,1.00){\makebox(0,0)[cc]{$_{3_0}$}}
\put(199.00,13.00){\makebox(0,0)[cc]{$_{\underline{3_1}}$}}
\put(177.00,16.00){\makebox(0,0)[cc]{$_{3_2}$}}
\put(168.00,28.00){\makebox(0,0)[cc]{$_{3_3}$}}
\put(165.00,45.00){\makebox(0,0)[cc]{$_{4_0}$}}
\put(165.00,62.00){\makebox(0,0)[cc]{$_{4_1}$}}
\put(169.00,80.00){\makebox(0,0)[cc]{$_{4_2}$}}
\put(177.00,92.00){\makebox(0,0)[cc]{$_{4_3}$}}
\emline{173.00}{30.00}{111}{194.00}{96.00}{112}
\emline{262.00}{30.00}{113}{242.00}{96.00}{114}
\emline{211.00}{6.00}{115}{265.00}{45.00}{116}
\emline{225.00}{6.00}{117}{171.00}{45.00}{118}
\put(82.00,54.00){\makebox(0,0)[cc]{$_{\mathcal P}$}}
\put(218.00,52.00){\makebox(0,0)[cc]{$_{\mathcal D}$}}
\put(259.00,92.00){\makebox(0,0)[cc]{$_{1_0}$}}
\put(247.00,7.00){\makebox(0,0)[cc]{$_{2_2}$}}
\put(230.00,1.00){\makebox(0,0)[cc]{$_{2_3}$}}
\put(194.00,97.00){\circle*{2.00}}
\put(226.00,101.00){\circle*{2.00}}
\put(254.00,89.00){\circle*{2.00}}
\put(182.00,17.00){\circle*{2.00}}
\put(262.00,77.00){\circle*{2.00}}
\put(266.00,45.00){\circle*{2.00}}
\put(242.00,97.00){\circle*{2.00}}
\put(194.00,10.00){\circle*{2.00}}
\end{picture}
\caption{Representations of ${\mathcal P}$ and ${\mathcal D}$}
\end{figure}
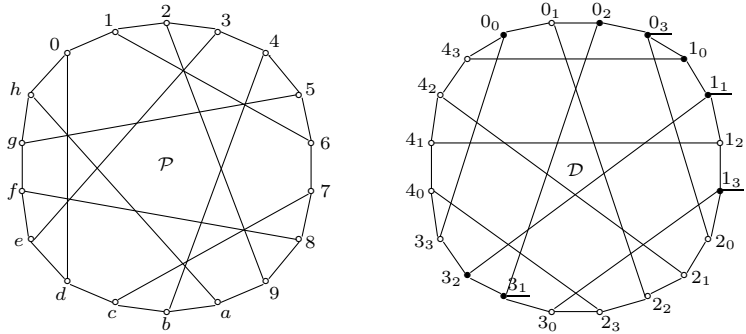

\proof For each positive integer $n$, let $I_n$ stand for the
$n$-cycle $(0,1,\ldots,n-1)$. ${\mathcal P}$ can be obtained from
$I_{18}$ by adding the edges $(1+6x,6+6x),(2+6x,9+6x),(4+6x,11+6x)$,
for $x\in\{0,1,2\}$, where operations are taken mod 18. Then $G$
admits the following collection of 6-cycles: $A_0=(123456)$,
$B_0=(3210de)$, $C_0=(34bcde)$, $D_0=(165gh0)$, $E_0=(329ab4)$,
(where octodecimal notation is used, up to $h=17$), as well as
$A_x,B_x,C_x,D_x,E_x$ obtained by uniformly adding $6x$ mod 18 to
the vertices of $A_0,B_0,C_0,D_0,E_0$, where
$x\in\Z_3\setminus\{0\}$, in addition to $F_0=(3298fe)$,
$F_1=(hg54ba)$, $F_2=(167cd0)$. These 6-cycles cannot be oriented
into a $(18\,\vec{C}_6)_{\vec{P}_3}$-OOC, for the following sequence
of alternating 6-cycles and 2-paths (with orientation reversed
between each 6-cycle and its corresponding succeeding 6-cycle)
reverses orientation from its initial 6-cycle to its terminal one:
$$\begin{array}{l}^{
D_1^{-1}654\,A_0123\,B_0210\,C_1h01\,D_0^{-1}g56\,C_2^{-1}765\,D_1}\\
^{=
(654bc7)\,654\,(123456)\,123\,(3210de)\,210\,(0129ah)\,h01\,(10hg56)\,g56\,
(5gf876)\,765(cb4567).}
\end{array}$$\noindent Another way to see this is via an auxiliary table for
$\mathcal P$, where $x=0,1,2$  (mod 3), presenting the form in which
the 6-cycles above share the 2-arcs, which are not always O-O for
$\mathcal P$, indicated by a minus sign in front of the heading of
each line of the table to distinguish it from the situation in
$\mathcal D$, shown below. Each $\eta_j$ in the table has subindex
$j$ indicating the equality of initial vertices $\eta_j=\xi_{i+2}$
of those 2-arcs, for $i=0,\ldots,5$:

\begin{eqnarray}\begin{array}{l|l}
^{-A_x:(\,B_x,E_x\,\,\,\,\,\,,E_{x+2},D_{x+1},D_x\,\,\,\,\,\,\,\,,B_{x+1})}
_{-B_x:(\,A_x,C_{x+1},F_2\,\,\,\,\,\,\,\,,A_{x+2},C_x\,\,\,\,\,\,\,\,,F_0\,\,\,\,\,\,\,\,)}
&^{-F_0:(E_0,B_1,E_1,B_2,E_2,B_0)}
_{-F_1:(D_0,E_2,D_1,E_0,D_2,E_1)}\\
^{-C_x:(\,E_x,D_{x+1},D_{x+2},B_{x+2},B_x\,\,\,\,\,\,\,\,,E_{x+2})}
_{-D_x:(A_x,C_{x+2}\,,F_1\,\,\,\,\,\,\,\,,A_{x+2},C_{x+1}\,\,,F_2\,\,\,\,\,\,\,)}
&^{-F_2:(B_1,D_1,B_2,D_2,B_0,D_0)}\\
^{-E_x:(F_0\,,C_{x+1}\,\,,A_{x+1},F_1\,\,\,\,\,\,\,\,,C_x\,\,\,\,\,\,\,\,,A_x\,\,\,\,\,\,)}&
\end{array}\end{eqnarray}
\noindent This proves that $\mathcal P$ is a
$\{\vec{C}_6\}^{\vec{P}_3}$-UH digraph.\bigskip

\noindent ${\mathcal D}$ can be obtained from $I_{20}$, with
vertices $4x,4x+1,4x+2,4x+3$ redenoted alternatively
$x_0,x_1,x_2,x_3$ respectively, for $x\in\Z_5$ by adding the edges
$(x_3,(x+2)_0)$ and $(x_1,(x+2)_2)$, with operations taken mod 5.
Then $G$ admits a $\{20\,\vec{C}_6\}_{\vec{P}_3}$-OOC formed by the
oriented 6-cycles $A_x,B_x,C_x,D_x$, for $x\in\{0,\ldots,4\}$, where
$$\begin{array}{l|l}
^{A_x=(x_0x_1x_2x_3(x+1)_0(x+4)_3)}_{C_x=(x_2x_1x_0(x+3)_3(x+3)_2(x+3)_1)} & ^{B_x=(x_1x_0(x+4)_3(x+4)_2(x+2)_1(x+2)_2)}_{D_x=(x_0(x+4)_3(x+1)_0(x+1)_1(x+3)_2(x+3)_3)}
\end{array}$$
The successive copies of $\vec{P}_3$ here, when reversed, in each
case, must belong to the following remaining oriented 6-cycles:
$$\begin{array}{l|l}
^{A_x:(C_x,C_{x+2},B_{x+1},D_{x+1},D_x,B_x)}_{C_x:(A_x,D_{x+4},D_x,A_{x+3},B_{x+1},B_{x+3})}&
^{B_x:(A_x,A_{x+4},D_{x+1},C_{x+4},C_{x+2},D_{x+4})}_{D_x:(A_x,C_{x+1},B_{x+1},B_{x+4},C_x,A_{x+4})}
\end{array}$$
showing that they constitute effectively an
$\{\eta\vec{C}_6\}_{\vec{P}_3}$-OOC.\qfd

\section{Cohering the distance-$2$ graphs of $6$-cycles}

\noindent We use the construction and notation of ${\mathcal D}$ and
its associated $\{\eta\vec{C}_6\}_{\vec{P}_3}$-OOC, as in the proof
of Theorem 2. Consider the collection $(\vec{\,\mathcal
C}_6)^2({\mathcal D})$ of squares of oriented $6$-cycles in the
$\{\eta\vec{C}_6\}_{\vec{P}_3}$-OOC of ${\mathcal D}$ in that proof.
Each arc $\vec{e}$ of a member $\vec{C}^2$ of $(\vec{\,\mathcal
C}_6)^2({\mathcal D})$ can be indicated by the middle vertex of the
2-arc $\vec{E}$ in $\vec{C}$ for which $\vec{e}$ stands, while the
tail and head of $\vec{e}$ are indicated by the tail and head of
$\vec{E}$, respectively. We cohere such $\vec{C}^2$s along their O-O
arc pairs in order to obtain a corresponding graph $Y({\mathcal D})$
with the $\{K_4,K_3\}_{K_2}$-UH property claimed in Subsection 1.4.
For such a setting, the following composed operation is performed,
where $\phi$ assigns to each 6-cycle in
$\{\eta\vec{C}_6\}_{\vec{P}_3}$-OOC its corresponding square:
\begin{equation}
{\mathcal
D}\,\,\,\rightarrow\,\,\,\{\eta\vec{C}_6\}_{\vec{P}_3}\mbox{-OOC}({\mathcal
D})\,\,\,\rightarrow^{\!\!\!\!\!\!\phi}\,\,\, (\vec{\,\mathcal
C}_6)^2({\mathcal D})\,\,\,\rightarrow\,\,\, Y({\mathcal
D}).\label{eqno1}\end{equation} We will explain in Section 4 how
this operation ${\mathcal D}\rightarrow Y({\mathcal D})$ is
performed.\bigskip

\noindent As mentioned in the table (1), in any oriented 6-cycle
$\xi$ of ${\mathcal P}$, each participating copy of $\vec{P}_3$,
when reversed in each case, must belong to a corresponding oriented
6-cycle $\eta$. In particular, each 6-cycle following such a copy of
$\vec{P}_3$ has its orientation reversed with respect to the one of
the preceding 6-cycle. This results in the second alternate 6-cycles
being considered with their orientation reversed with respect to the
first alternate 6-cycles. Because of this, we say that there are $2$
{\it alternate} O-O $\{\frac{1}{2}\eta\vec{C}_6\}_{\vec{P}_3}$-OOCs,
in the absence of just one $\{\eta\vec{C}_6\}_{\vec{P}_3}$-OOC for
$\mathcal P$. This allows $2$ corresponding alternate
half-operations similar in nature to (2), above. See Section 6
below.\bigskip

\noindent The $2$ versions of $\vec{\,\mathcal C}_6^{2}({\mathcal
P})$ here and the only one of $Y({\mathcal D})=\vec{\,\mathcal
B}_6^{2}({\mathcal D})$ are formed by oriented triangles that
determine $2$ corresponding graphs $Y_1({\mathcal P})$ and
$Y_2({\mathcal P})$ and a single graph $Y({\mathcal D})$ with $2$
components $Y_1({\mathcal D})$ and $Y_2({\mathcal D})$.

\section{Desargues reattachment Menger graph}

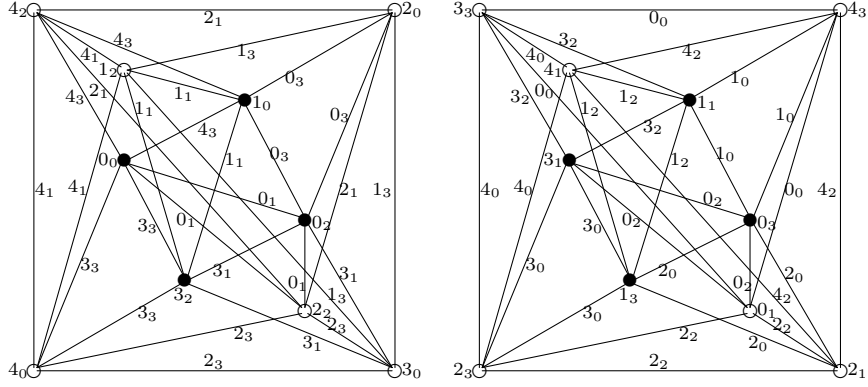
\begin{figure}[htp]
\unitlength=0.4mm \special{em:linewidth 0.4pt} \linethickness{0.4pt}
\begin{picture}(280.00,124.00)
\emline{7.00}{2.00}{1}{125.00}{2.00}{2}
\emline{126.00}{3.00}{3}{126.00}{121.00}{4}
\emline{125.00}{122.00}{5}{7.00}{122.00}{6}
\emline{6.00}{121.00}{7}{6.00}{3.00}{8}
\put(132.00,2.00){\makebox(0,0)[cc]{$_{3_0}$}}
\put(1.00,122.00){\makebox(0,0)[cc]{$_{4_2}$}}
\emline{37.00}{71.00}{9}{95.00}{53.00}{10}
\emline{96.00}{53.00}{11}{76.00}{91.00}{12}
\emline{75.00}{91.00}{13}{57.00}{33.00}{14}
\emline{75.00}{92.00}{15}{37.00}{72.00}{16}
\emline{36.00}{71.00}{17}{56.00}{33.00}{18}
\emline{57.00}{32.00}{19}{95.00}{51.00}{20}
\emline{97.00}{53.00}{21}{125.00}{121.00}{22}
\emline{97.00}{51.00}{23}{125.00}{3.00}{24}
\emline{125.00}{121.00}{25}{77.00}{93.00}{26}
\emline{7.00}{3.00}{27}{55.00}{31.00}{28}
\emline{57.00}{31.00}{29}{125.00}{3.00}{30}
\emline{35.00}{71.00}{31}{7.00}{3.00}{32}
\emline{7.00}{121.00}{33}{35.00}{73.00}{34}
\put(1.00,2.00){\makebox(0,0)[cc]{$_{4_0}$}}
\put(132.00,122.00){\makebox(0,0)[cc]{$_{2_0}$}}
\put(31.00,72.00){\makebox(0,0)[cc]{$_{0_0}$}}
\put(82.00,91.00){\makebox(0,0)[cc]{$_{1_0}$}}
\put(56.00,27.00){\makebox(0,0)[cc]{$_{3_2}$}}
\put(102.00,51.00){\makebox(0,0)[cc]{$_{0_2}$}}
\put(102.00,22.00){\makebox(0,0)[cc]{$_{2_2}$}}
\put(31.00,102.00){\makebox(0,0)[cc]{$_{1_2}$}}
\emline{96.00}{51.00}{35}{96.00}{23.00}{36}
\emline{96.00}{23.00}{37}{125.00}{121.00}{38}
\emline{125.00}{3.00}{39}{97.00}{21.00}{40}
\emline{95.00}{21.00}{41}{7.00}{3.00}{42}
\emline{35.00}{103.00}{43}{7.00}{121.00}{44}
\emline{95.00}{23.00}{45}{37.00}{71.00}{46}
\emline{37.00}{102.00}{47}{75.00}{92.00}{48}
\emline{36.00}{101.00}{49}{56.00}{33.00}{50}
\emline{95.00}{23.00}{51}{7.00}{121.00}{52}
\emline{37.00}{101.00}{53}{125.00}{3.00}{54}
\put(66.00,119.00){\makebox(0,0)[cc]{$_{2_1}$}}
\put(122.00,62.00){\makebox(0,0)[cc]{$_{1_3}$}}
\put(66.00,5.00){\makebox(0,0)[cc]{$_{2_3}$}}
\put(10.00,62.00){\makebox(0,0)[cc]{$_{4_1}$}}
\put(20.00,93.00){\makebox(0,0)[cc]{$_{4_3}$}}
\put(28.00,95.00){\makebox(0,0)[cc]{$_{2_1}$}}
\put(25.00,107.00){\makebox(0,0)[cc]{$_{4_1}$}}
\put(36.00,112.00){\makebox(0,0)[cc]{$_{4_3}$}}
\put(77.00,108.00){\makebox(0,0)[cc]{$_{1_3}$}}
\put(93.00,99.00){\makebox(0,0)[cc]{$_{0_3}$}}
\put(108.00,87.00){\makebox(0,0)[cc]{$_{0_3}$}}
\put(88.00,74.00){\makebox(0,0)[cc]{$_{0_3}$}}
\put(111.00,62.00){\makebox(0,0)[cc]{$_{2_1}$}}
\put(111.00,34.00){\makebox(0,0)[cc]{$_{3_1}$}}
\put(107.00,27.00){\makebox(0,0)[cc]{$_{1_3}$}}
\put(107.00,17.00){\makebox(0,0)[cc]{$_{2_3}$}}
\put(99.00,11.00){\makebox(0,0)[cc]{$_{3_1}$}}
\put(76.00,14.00){\makebox(0,0)[cc]{$_{2_3}$}}
\put(44.00,21.00){\makebox(0,0)[cc]{$_{3_3}$}}
\put(25.00,37.00){\makebox(0,0)[cc]{$_{3_3}$}}
\put(73.00,72.00){\makebox(0,0)[cc]{$_{1_1}$}}
\put(56.00,94.00){\makebox(0,0)[cc]{$_{1_1}$}}
\put(43.00,89.00){\makebox(0,0)[cc]{$_{1_1}$}}
\put(84.00,59.00){\makebox(0,0)[cc]{$_{0_1}$}}
\put(94.00,31.00){\makebox(0,0)[cc]{$_{0_1}$}}
\put(57.00,52.00){\makebox(0,0)[cc]{$_{0_1}$}}
\put(64.00,82.00){\makebox(0,0)[cc]{$_{4_3}$}}
\put(44.00,50.00){\makebox(0,0)[cc]{$_{3_3}$}}
\put(69.00,35.00){\makebox(0,0)[cc]{$_{3_1}$}}
\emline{7.00}{121.00}{55}{75.00}{93.00}{56}
\emline{155.00}{2.00}{57}{273.00}{2.00}{58}
\emline{274.00}{3.00}{59}{274.00}{121.00}{60}
\emline{273.00}{122.00}{61}{155.00}{122.00}{62}
\emline{154.00}{121.00}{63}{154.00}{3.00}{64}
\put(280.00,2.00){\makebox(0,0)[cc]{$_{2_1}$}}
\put(149.00,122.00){\makebox(0,0)[cc]{$_{3_3}$}}
\emline{185.00}{71.00}{65}{243.00}{53.00}{66}
\emline{244.00}{53.00}{67}{224.00}{91.00}{68}
\emline{223.00}{91.00}{69}{205.00}{33.00}{70}
\emline{223.00}{92.00}{71}{185.00}{72.00}{72}
\emline{184.00}{71.00}{73}{204.00}{33.00}{74}
\emline{205.00}{32.00}{75}{243.00}{51.00}{76}
\emline{245.00}{53.00}{77}{273.00}{121.00}{78}
\emline{245.00}{51.00}{79}{273.00}{3.00}{80}
\emline{273.00}{121.00}{81}{225.00}{93.00}{82}
\emline{155.00}{3.00}{83}{203.00}{31.00}{84}
\emline{205.00}{31.00}{85}{273.00}{3.00}{86}
\emline{183.00}{71.00}{87}{155.00}{3.00}{88}
\emline{155.00}{121.00}{89}{183.00}{73.00}{90}
\put(149.00,2.00){\makebox(0,0)[cc]{$_{2_3}$}}
\put(280.00,122.00){\makebox(0,0)[cc]{$_{4_3}$}}
\put(179.00,72.00){\makebox(0,0)[cc]{$_{3_1}$}}
\put(230.00,91.00){\makebox(0,0)[cc]{$_{1_1}$}}
\put(204.00,27.00){\makebox(0,0)[cc]{$_{1_3}$}}
\put(250.00,51.00){\makebox(0,0)[cc]{$_{0_3}$}}
\put(250.00,22.00){\makebox(0,0)[cc]{$_{0_1}$}}
\put(179.00,102.00){\makebox(0,0)[cc]{$_{4_1}$}}
\emline{244.00}{51.00}{91}{244.00}{23.00}{92}
\emline{244.00}{23.00}{93}{273.00}{121.00}{94}
\emline{273.00}{3.00}{95}{245.00}{21.00}{96}
\emline{243.00}{21.00}{97}{155.00}{3.00}{98}
\emline{183.00}{103.00}{99}{155.00}{121.00}{100}
\emline{243.00}{23.00}{101}{185.00}{71.00}{102}
\emline{185.00}{102.00}{103}{223.00}{92.00}{104}
\emline{184.00}{101.00}{105}{204.00}{33.00}{106}
\emline{243.00}{23.00}{107}{155.00}{121.00}{108}
\emline{185.00}{101.00}{109}{273.00}{3.00}{110}
\put(214.00,119.00){\makebox(0,0)[cc]{$_{0_0}$}}
\put(270.00,62.00){\makebox(0,0)[cc]{$_{4_2}$}}
\put(214.00,5.00){\makebox(0,0)[cc]{$_{2_2}$}}
\put(158.00,62.00){\makebox(0,0)[cc]{$_{4_0}$}}
\put(168.00,93.00){\makebox(0,0)[cc]{$_{3_2}$}}
\put(176.00,95.00){\makebox(0,0)[cc]{$_{0_0}$}}
\put(173.00,107.00){\makebox(0,0)[cc]{$_{4_0}$}}
\put(184.00,112.00){\makebox(0,0)[cc]{$_{3_2}$}}
\put(225.00,108.00){\makebox(0,0)[cc]{$_{4_2}$}}
\put(241.00,99.00){\makebox(0,0)[cc]{$_{1_0}$}}
\put(256.00,87.00){\makebox(0,0)[cc]{$_{1_0}$}}
\put(236.00,74.00){\makebox(0,0)[cc]{$_{1_0}$}}
\put(259.00,62.00){\makebox(0,0)[cc]{$_{0_0}$}}
\put(259.00,34.00){\makebox(0,0)[cc]{$_{2_0}$}}
\put(255.00,27.00){\makebox(0,0)[cc]{$_{4_2}$}}
\put(255.00,17.00){\makebox(0,0)[cc]{$_{2_2}$}}
\put(247.00,11.00){\makebox(0,0)[cc]{$_{2_0}$}}
\put(224.00,14.00){\makebox(0,0)[cc]{$_{2_2}$}}
\put(192.00,21.00){\makebox(0,0)[cc]{$_{3_0}$}}
\put(173.00,37.00){\makebox(0,0)[cc]{$_{3_0}$}}
\put(221.00,72.00){\makebox(0,0)[cc]{$_{1_2}$}}
\put(204.00,93.00){\makebox(0,0)[cc]{$_{1_2}$}}
\put(191.00,89.00){\makebox(0,0)[cc]{$_{1_2}$}}
\put(232.00,59.00){\makebox(0,0)[cc]{$_{0_2}$}}
\put(242.00,31.00){\makebox(0,0)[cc]{$_{0_2}$}}
\put(205.00,52.00){\makebox(0,0)[cc]{$_{0_2}$}}
\put(212.00,83.00){\makebox(0,0)[cc]{$_{3_2}$}}
\put(192.00,50.00){\makebox(0,0)[cc]{$_{3_0}$}}
\put(217.00,35.00){\makebox(0,0)[cc]{$_{2_0}$}}
\emline{155.00}{121.00}{111}{223.00}{93.00}{112}
\emline{125.00}{121.00}{113}{37.00}{102.00}{114}
\emline{273.00}{121.00}{115}{185.00}{102.00}{116}
\emline{183.00}{101.00}{117}{155.00}{3.00}{118}
\put(169.00,63.00){\makebox(0,0)[cc]{$_{4_0}$}}
\emline{35.00}{101.00}{119}{7.00}{3.00}{120}
\put(21.00,63.00){\makebox(0,0)[cc]{$_{4_1}$}}
\put(36.00,72.00){\circle*{4.00}}
\put(96.00,52.00){\circle*{4.00}}
\put(56.00,32.00){\circle*{4.00}}
\put(76.00,92.00){\circle*{4.00}}
\put(6.00,122.00){\circle{4.00}}
\put(126.00,122.00){\circle{4.00}}
\put(126.00,2.00){\circle{4.00}}
\put(6.00,2.00){\circle{4.00}}
\put(36.00,102.00){\circle{4.00}}
\put(96.00,22.00){\circle{4.00}}
\put(184.00,72.00){\circle*{4.00}}
\put(224.00,92.00){\circle*{4.00}}
\put(204.00,32.00){\circle*{4.00}}
\put(244.00,52.00){\circle*{4.00}}
\put(154.00,2.00){\circle{4.00}}
\put(154.00,122.00){\circle{4.00}}
\put(274.00,122.00){\circle{4.00}}
\put(274.00,2.00){\circle{4.00}}
\put(244.00,22.00){\circle{4.00}}
\put(184.00,102.00){\circle{4.00}}
\end{picture}
\caption{Representations of $Y_1({\mathcal D})$ and $Y_2({\mathcal D})$}
\end{figure}

\noindent For $i=1,2$, it is a matter of checking that
$Y_i({\mathcal D})$ is an isomorphic coherent $\{K_4,K_3\}_{K_2}$-UH
graph formed by 5 copies of $K_4$ and 10 of $K_3\not\subset K_4$,
with each such copy of $K_3$
%: {\bf(a)} not forming part of a copy of $K_4$ in $Y_i(G)$; {\bf(b)}
having its edges indicated by a constant symbol, as shown in Figure
2. Each of the 5 copies $\mathcal T$ of $K_4$ in $Y_i({\mathcal D})$
has any one of its six edges as a pair of O-O arcs, say $\vec{e}$
and $(\vec{e})^{-1}$, arising from corresponding O-O 2-arcs
separating $2$ oriented 6-cycles of the
$\{20\,\vec{C}_6\}_{\vec{P}_3}\mbox{-OOC}({\mathcal D})$, as
obtained in the proof of Theorem 2. Moreover, these $2$ oriented
6-cycles have images in $(\vec{\,\mathcal C}_6)^2(G)$ via the map
$\phi$ displayed in (2) that are the $2$ oriented triangles that
share $e$ in $\mathcal T$, or in other words, having $\vec{e}$ and
$(\vec{e})^{-1}$ as O-O arcs.\bigskip

\noindent For example, each of the $2$ `central' copies of $K_4$
with black vertices in either side of Figure 2 has the 4 composing
vertex triples, namely $(0_0,0_2,1_0)$, $(0_2,1_0,3_2)$,
$(1_0,3_2,0_0)$ and $(3_2,0_0,0_2)$, for $Y_1({\mathcal D})$, (resp.
$(0_3,1_1,1_3)$,

\noindent $(1_1,1_3,3_1)$, $(1_3,3_1,0_3)$ and $(3_1,0_3,1_1)$, for
$Y_2({\mathcal D})$) as alternate vertices of 4 corresponding
6-cycles in ${\mathcal D}$, as can be checked on the right side of
Figure 1, where the corresponding vertices in $\mathcal D$ are also
black and those corresponding to vertices of $Y_2({\mathcal D})$
underlined for distinction. Now, the edge $(0_0,0_2)$ in
$Y_1({\mathcal D})$, corresponding to the 2-path $(0_0,0_1,0_2)$ in
${\mathcal D}$, has its $2$ composing arcs separating the oriented
triangles $\phi(A_0)=(0_0,0_2,1_0)$ and $\phi(C_0)=(0_0,0_2,3_2)$,
corresponding to the oriented 6-cycles
$A_0=(0_0,0_1,0_2,0_3,1_0,4_3)$ and
$C_0=(0_2,0_1,0_0,3_3,3_2,3_1)$.\bigskip

\noindent We notice that the 10 vertices and 10 copies of
$K_3\not\subset K_4$ in either $Y_i({\mathcal D})$, ($i=1,2$), may
be considered as the points and lines of the Desargues self-dual
$(10_3)$ configuration, and that the Menger graph of this coincides
with $Y_i({\mathcal D})$ \cite{Cox}. Each vertex of $Y_i({\mathcal
D})$ is the meeting vertex of $2$ copies of $K_4$ and 3 copies of
$K_3$ not forming part of a copy of $K_4$.

\begin{thm}
$Y_1({\mathcal D})$ and $Y_2({\mathcal D})$ are coherent
$\{K_4,K_3\}_{K_2}$-UH graphs composed by $5$ copies of $K_4$ and
$10$ copies of $K_3\not\subset K_4$ each. Moreover, the $10$
vertices and $10$ copies of $K_3\not\subset K_4$ in either graph
constitute the Desargues self-dual $(10_3)$ con\-fi\-gu\-ra\-tion
whose Levi graph is $\mathcal D$ and whose Menger graph is equal to both $Y_1({\mathcal D})$ and $Y_2({\mathcal D})$.
Furthermore, both
graphs are isomorphic to $L(K_5)$, whose complement is the Petersen
graph. \qfd\end{thm}

\noindent Deleting a copy $H$ of $K_4$ from such $Y_i({\mathcal D})$
($i=1,2$)  yields a copy $J$ of $K_{2,2,2}$\,, 4 of whose composing
copies of $K_3$, with no common edges, are faces of corres\-ponding
copies of $K_4\neq H$. The other 4 copies of $K_3$ are among the 10
mentioned copies of $K_3$ in $G$. A realization of $Y_i({\mathcal
D})$ in 3-space can be obtained from a regular octahedron $O_3$ with
1-skeleton $J$ via the midpoints, say $x_1,x_2,x_3,x_4$, of the 4
segments joining the barycenters of 4 edge-disjoint alternate
triangles, say $T_1,T_2,T_3,T_4$, in $O_3$ to the barycenter of
$O_3$: just construct the tetrahedron $\Delta_j$ determined by each
$T_j$ and corresponding $x_j$, as well as the tetrahedron $\Delta_0$
determined by the 4 $x_i$\thinspace s.\bigskip

\noindent A realization $\kappa$ of $K_5$ in 3-space is obtained
whose vertices are the barycenters of
$\Delta_0,\Delta_1,\Delta_2,\Delta_3,\Delta_4$ and whose edges are
the segments that join those barycenters. By taking the midpoints of
the segments realizing the edges of $\kappa$ and joining each two of
them, say midpoints $P$ and $Q$ of respective segments $p$ and $q$,
by a new segment whenever $p$ and $q$ have an end in common in
$\kappa$, a realization $L(\kappa)$ of $L(K_5)$ is obtained. This
$L(\kappa)$ is a smaller realization of $L(K_5)$ than that of
$Y_i({\mathcal D})$ in the previous paragraph and leads to an
octahedron $O'_3\subset O_3$ by the deletion of its central copy of
$K_4$. This procedure may be repeated indefinitely, generating a
nested sequence of realizations of $Y_i({\mathcal D})$ in 3-space.
Since $Y_1({\mathcal D})$ and $Y_2({\mathcal D})$ are isomorphic to
$L(K_5)$, whose complement is the Petersen graph, this sequence
yields a corresponding infinite sequence of realizations of the
Petersen graph in 3-space.

\section{Generalization of Theorem 3}

\noindent Theorem 3 can be partly generalized by replacing $L(K_5)$
by $L(K_n)$ ($n\geq 4$). This produces a coherent
$\{K_{n-1},K_3\}_{K_2}$-UH graph.

\begin{thm}
The line graph $L(K_n)$, with $n\geq 4$, is a coherent $\{K_{n-1},$
$K_3\}_{K_2}$-UH graph with $n$ copies of $K_{n-1}$ and ${n\choose
3}$ copies of $K_3\not\subset K_{n-1}$.\end{thm}

\proof Each vertex $v$ of $K_n$ is taken as a color of edges of
$L(K_n)$ under the following rule: color all the edges between
vertices of $L(K_n)$ representing edges incident to $v$ with color
$v$. Then, each triple of edge colors of $L(K_n)$ corresponds to the
edges of a well determined copy of $K_3\not\subset K_{n-1}$ in
$L(K_n)$. Thus, there are exactly ${n\choose 3}$ copies of
$K_3\not\subset K_{n-1}$ intervening in $L(K_n)$ looked upon as a
coherent $\{K_{n-1},K_3\}_{K_2}$-UH graph.\qfd

\section{Pappus reattachment Menger graph}

\begin{figure}[htp]
\unitlength=0.4mm \special{em:inewidth 0.4pt} \linethickness{0.4pt}
\begin{picture}(280.00,122.00)
\put(5.00,1.00){\circle{2.00}}
\put(125.00,1.00){\circle{2.00}}
\put(5.00,121.00){\circle{2.00}}
\put(125.00,121.00){\circle{2.00}}
\put(130.00,1.00){\makebox(0,0)[cc]{$_b$}}
\put(1.00,121.00){\makebox(0,0)[cc]{$_b$}}
\put(1.00,1.00){\makebox(0,0)[cc]{$_b$}}
\put(130.00,121.00){\makebox(0,0)[cc]{$_b$}}
\put(45.00,1.00){\circle{2.00}}
\put(85.00,1.00){\circle{2.00}}
\put(45.00,121.00){\circle{2.00}}
\put(85.00,121.00){\circle{2.00}}
\put(5.00,41.00){\circle{2.00}}
\put(125.00,41.00){\circle{2.00}}
\put(45.00,41.00){\circle{2.00}}
\put(85.00,41.00){\circle{2.00}}
\put(5.00,81.00){\circle{2.00}}
\put(125.00,81.00){\circle{2.00}}
\put(45.00,81.00){\circle{2.00}}
\put(85.00,81.00){\circle{2.00}}
\emline{6.00}{1.00}{1}{44.00}{1.00}{2}
\emline{6.00}{41.00}{3}{44.00}{41.00}{4}
\emline{6.00}{81.00}{5}{44.00}{81.00}{6}
\emline{6.00}{121.00}{7}{44.00}{121.00}{8}
\emline{44.00}{2.00}{9}{6.00}{40.00}{10}
\emline{5.00}{2.00}{11}{5.00}{40.00}{12}
\emline{44.00}{42.00}{13}{6.00}{80.00}{14}
\emline{5.00}{42.00}{15}{5.00}{80.00}{16}
\emline{44.00}{82.00}{17}{6.00}{120.00}{18}
\emline{5.00}{82.00}{19}{5.00}{120.00}{20}
\emline{46.00}{1.00}{21}{84.00}{1.00}{22}
\emline{46.00}{41.00}{23}{84.00}{41.00}{24}
\emline{46.00}{81.00}{25}{84.00}{81.00}{26}
\emline{46.00}{121.00}{27}{84.00}{121.00}{28}
\emline{84.00}{2.00}{29}{46.00}{40.00}{30}
\emline{45.00}{2.00}{31}{45.00}{40.00}{32}
\emline{84.00}{42.00}{33}{46.00}{80.00}{34}
\emline{45.00}{42.00}{35}{45.00}{80.00}{36}
\emline{84.00}{82.00}{37}{46.00}{120.00}{38}
\emline{45.00}{82.00}{39}{45.00}{120.00}{40}
\emline{86.00}{1.00}{41}{124.00}{1.00}{42}
\emline{86.00}{41.00}{43}{124.00}{41.00}{44}
\emline{86.00}{81.00}{45}{124.00}{81.00}{46}
\emline{86.00}{121.00}{47}{124.00}{121.00}{48}
\emline{124.00}{2.00}{49}{86.00}{40.00}{50}
\emline{85.00}{2.00}{51}{85.00}{40.00}{52}
\emline{124.00}{42.00}{53}{86.00}{80.00}{54}
\emline{85.00}{42.00}{55}{85.00}{80.00}{56}
\emline{124.00}{82.00}{57}{86.00}{120.00}{58}
\emline{85.00}{82.00}{59}{85.00}{120.00}{60}
\emline{125.00}{2.00}{61}{125.00}{40.00}{62}
\emline{125.00}{42.00}{63}{125.00}{80.00}{64}
\emline{125.00}{82.00}{65}{125.00}{120.00}{66}
\put(130.00,41.00){\makebox(0,0)[cc]{$_d$}}
\put(1.00,41.00){\makebox(0,0)[cc]{$_d$}}
\put(130.00,81.00){\makebox(0,0)[cc]{$_7$}}
\put(1.00,81.00){\makebox(0,0)[cc]{$_7$}}
\put(42.00,118.00){\makebox(0,0)[cc]{$_3$}}
\put(42.00,-2.00){\makebox(0,0)[cc]{$_3$}}
\put(42.00,38.00){\makebox(0,0)[cc]{$_1$}}
\put(42.00,78.00){\makebox(0,0)[cc]{$_9$}}
\put(82.00,118.00){\makebox(0,0)[cc]{$_5$}}
\put(82.00,-2.00){\makebox(0,0)[cc]{$_5$}}
\put(82.00,38.00){\makebox(0,0)[cc]{$_h$}}
\put(82.00,78.00){\makebox(0,0)[cc]{$_f$}}
\put(1.00,21.00){\makebox(0,0)[cc]{$_c$}}
\put(1.00,61.00){\makebox(0,0)[cc]{$_c$}}
\put(1.00,101.00){\makebox(0,0)[cc]{$_c$}}
\put(130.00,21.00){\makebox(0,0)[cc]{$_c$}}
\put(130.00,61.00){\makebox(0,0)[cc]{$_c$}}
\put(130.00,101.00){\makebox(0,0)[cc]{$_c$}}
\put(41.00,21.00){\makebox(0,0)[cc]{$_2$}}
\put(41.00,61.00){\makebox(0,0)[cc]{$_2$}}
\put(41.00,101.00){\makebox(0,0)[cc]{$_2$}}
\put(81.00,21.00){\makebox(0,0)[cc]{$_g$}}
\put(81.00,61.00){\makebox(0,0)[cc]{$_g$}}
\put(81.00,101.00){\makebox(0,0)[cc]{$_g$}}
\put(21.00,21.00){\makebox(0,0)[cc]{$_e$}}
\put(21.00,61.00){\makebox(0,0)[cc]{$_6$}}
\put(21.00,101.00){\makebox(0,0)[cc]{$_a$}}
\put(61.00,21.00){\makebox(0,0)[cc]{$_6$}}
\put(61.00,61.00){\makebox(0,0)[cc]{$_a$}}
\put(61.00,101.00){\makebox(0,0)[cc]{$_e$}}
\put(101.00,21.00){\makebox(0,0)[cc]{$_a$}}
\put(101.00,61.00){\makebox(0,0)[cc]{$_e$}}
\put(101.00,101.00){\makebox(0,0)[cc]{$_6$}}
\put(25.00,-2.00){\makebox(0,0)[cc]{$_4$}}
\put(65.00,-2.00){\makebox(0,0)[cc]{$_4$}}
\put(105.00,-2.00){\makebox(0,0)[cc]{$_4$}}
\put(25.00,38.00){\makebox(0,0)[cc]{$_0$}}
\put(65.00,38.00){\makebox(0,0)[cc]{$_0$}}
\put(105.00,38.00){\makebox(0,0)[cc]{$_0$}}
\put(25.00,78.00){\makebox(0,0)[cc]{$_8$}}
\put(65.00,78.00){\makebox(0,0)[cc]{$_8$}}
\put(105.00,78.00){\makebox(0,0)[cc]{$_8$}}
\put(25.00,118.00){\makebox(0,0)[cc]{$_4$}}
\put(65.00,118.00){\makebox(0,0)[cc]{$_4$}}
\put(105.00,118.00){\makebox(0,0)[cc]{$_4$}}
\put(30.00,26.00){\makebox(0,0)[cc]{$_{B_0^1}$}}
\put(70.00,26.00){\makebox(0,0)[cc]{$_{D_0^1}$}}
\put(110.00,26.00){\makebox(0,0)[cc]{$_{D_2^1}$}}
\put(30.00,66.00){\makebox(0,0)[cc]{$_{B_1^1}$}}
\put(70.00,66.00){\makebox(0,0)[cc]{$_{E_1^1}$}}
\put(110.00,66.00){\makebox(0,0)[cc]{$_{B_2^1}$}}
\put(30.00,106.00){\makebox(0,0)[cc]{$_{E_0^1}$}}
\put(70.00,106.00){\makebox(0,0)[cc]{$_{E_2^1}$}}
\put(110.00,106.00){\makebox(0,0)[cc]{$_{D_1^1}$}}
\put(15.00,11.00){\makebox(0,0)[cc]{$_{C_0^1}$}}
\put(55.00,11.00){\makebox(0,0)[cc]{$_{A_0^1}$}}
\put(95.00,11.00){\makebox(0,0)[cc]{$_{F_1^1}$}}
\put(15.00,51.00){\makebox(0,0)[cc]{$_{F_2^1}$}}
\put(55.00,51.00){\makebox(0,0)[cc]{$_{C_1^1}$}}
\put(95.00,51.00){\makebox(0,0)[cc]{$_{A_2^1}$}}
\put(15.00,91.00){\makebox(0,0)[cc]{$_{A_1^1}$}}
\put(55.00,91.00){\makebox(0,0)[cc]{$_{F_0^1}$}}
\put(95.00,91.00){\makebox(0,0)[cc]{$_{C_2^1}$}}
\put(155.00,1.00){\circle{2.00}}
\put(275.00,1.00){\circle{2.00}}
\put(155.00,121.00){\circle{2.00}}
\put(275.00,121.00){\circle{2.00}}
\put(280.00,1.00){\makebox(0,0)[cc]{$_e$}}
\put(151.00,121.00){\makebox(0,0)[cc]{$_e$}}
\put(151.00,1.00){\makebox(0,0)[cc]{$_e$}}
\put(280.00,121.00){\makebox(0,0)[cc]{$_e$}}
\put(195.00,1.00){\circle{2.00}}
\put(235.00,1.00){\circle{2.00}}
\put(195.00,121.00){\circle{2.00}}
\put(235.00,121.00){\circle{2.00}}
\put(155.00,41.00){\circle{2.00}}
\put(275.00,41.00){\circle{2.00}}
\put(195.00,41.00){\circle{2.00}}
\put(235.00,41.00){\circle{2.00}}
\put(155.00,81.00){\circle{2.00}}
\put(275.00,81.00){\circle{2.00}}
\put(195.00,81.00){\circle{2.00}}
\put(235.00,81.00){\circle{2.00}}
\emline{156.00}{1.00}{67}{194.00}{1.00}{68}
\emline{156.00}{41.00}{69}{194.00}{41.00}{70}
\emline{156.00}{81.00}{71}{194.00}{81.00}{72}
\emline{156.00}{121.00}{73}{194.00}{121.00}{74}
\emline{194.00}{2.00}{75}{156.00}{40.00}{76}
\emline{155.00}{2.00}{77}{155.00}{40.00}{78}
\emline{194.00}{42.00}{79}{156.00}{80.00}{80}
\emline{155.00}{42.00}{81}{155.00}{80.00}{82}
\emline{194.00}{82.00}{83}{156.00}{120.00}{84}
\emline{155.00}{82.00}{85}{155.00}{120.00}{86}
\emline{196.00}{1.00}{87}{234.00}{1.00}{88}
\emline{196.00}{41.00}{89}{234.00}{41.00}{90}
\emline{196.00}{81.00}{91}{234.00}{81.00}{92}
\emline{196.00}{121.00}{93}{234.00}{121.00}{94}
\emline{234.00}{2.00}{95}{196.00}{40.00}{96}
\emline{195.00}{2.00}{97}{195.00}{40.00}{98}
\emline{234.00}{42.00}{99}{196.00}{80.00}{100}
\emline{195.00}{42.00}{101}{195.00}{80.00}{102}
\emline{234.00}{82.00}{103}{196.00}{120.00}{104}
\emline{195.00}{82.00}{105}{195.00}{120.00}{106}
\emline{236.00}{1.00}{107}{274.00}{1.00}{108}
\emline{236.00}{41.00}{109}{274.00}{41.00}{110}
\emline{236.00}{81.00}{111}{274.00}{81.00}{112}
\emline{236.00}{121.00}{113}{274.00}{121.00}{114}
\emline{274.00}{2.00}{115}{236.00}{40.00}{116}
\emline{235.00}{2.00}{117}{235.00}{40.00}{118}
\emline{274.00}{42.00}{119}{236.00}{80.00}{120}
\emline{235.00}{42.00}{121}{235.00}{80.00}{122}
\emline{274.00}{82.00}{123}{236.00}{120.00}{124}
\emline{235.00}{82.00}{125}{235.00}{120.00}{126}
\emline{275.00}{2.00}{127}{275.00}{40.00}{128}
\emline{275.00}{42.00}{129}{275.00}{80.00}{130}
\emline{275.00}{82.00}{131}{275.00}{120.00}{132}
\put(280.00,41.00){\makebox(0,0)[cc]{$_g$}}
\put(151.00,41.00){\makebox(0,0)[cc]{$_g$}}
\put(280.00,81.00){\makebox(0,0)[cc]{$_8$}}
\put(151.00,81.00){\makebox(0,0)[cc]{$_8$}}
\put(192.00,118.00){\makebox(0,0)[cc]{$_4$}}
\put(192.00,-2.00){\makebox(0,0)[cc]{$_4$}}
\put(192.00,38.00){\makebox(0,0)[cc]{$_a$}}
\put(192.00,78.00){\makebox(0,0)[cc]{$_c$}}
\put(232.00,118.00){\makebox(0,0)[cc]{$_2$}}
\put(232.00,-2.00){\makebox(0,0)[cc]{$_2$}}
\put(232.00,38.00){\makebox(0,0)[cc]{$_0$}}
\put(232.00,78.00){\makebox(0,0)[cc]{$_6$}}
\put(151.00,21.00){\makebox(0,0)[cc]{$_f$}}
\put(151.00,61.00){\makebox(0,0)[cc]{$_f$}}
\put(151.00,101.00){\makebox(0,0)[cc]{$_f$}}
\put(280.00,21.00){\makebox(0,0)[cc]{$_f$}}
\put(280.00,61.00){\makebox(0,0)[cc]{$_f$}}
\put(280.00,101.00){\makebox(0,0)[cc]{$_f$}}
\put(191.00,21.00){\makebox(0,0)[cc]{$_b$}}
\put(191.00,61.00){\makebox(0,0)[cc]{$_b$}}
\put(191.00,101.00){\makebox(0,0)[cc]{$_b$}}
\put(231.00,21.00){\makebox(0,0)[cc]{$_1$}}
\put(231.00,61.00){\makebox(0,0)[cc]{$_1$}}
\put(231.00,101.00){\makebox(0,0)[cc]{$_1$}}
\put(171.00,21.00){\makebox(0,0)[cc]{$_5$}}
\put(171.00,61.00){\makebox(0,0)[cc]{$_9$}}
\put(171.00,101.00){\makebox(0,0)[cc]{$_d$}}
\put(211.00,21.00){\makebox(0,0)[cc]{$_9$}}
\put(211.00,61.00){\makebox(0,0)[cc]{$_d$}}
\put(211.00,101.00){\makebox(0,0)[cc]{$_5$}}
\put(251.00,21.00){\makebox(0,0)[cc]{$_d$}}
\put(251.00,61.00){\makebox(0,0)[cc]{$_5$}}
\put(251.00,101.00){\makebox(0,0)[cc]{$_9$}}
\put(175.00,-2.00){\makebox(0,0)[cc]{$_3$}}
\put(215.00,-2.00){\makebox(0,0)[cc]{$_3$}}
\put(255.00,-2.00){\makebox(0,0)[cc]{$_3$}}
\put(175.00,38.00){\makebox(0,0)[cc]{$_h$}}
\put(215.00,38.00){\makebox(0,0)[cc]{$_h$}}
\put(255.00,38.00){\makebox(0,0)[cc]{$_h$}}
\put(175.00,78.00){\makebox(0,0)[cc]{$_7$}}
\put(215.00,78.00){\makebox(0,0)[cc]{$_7$}}
\put(255.00,78.00){\makebox(0,0)[cc]{$_7$}}
\put(175.00,118.00){\makebox(0,0)[cc]{$_3$}}
\put(215.00,118.00){\makebox(0,0)[cc]{$_3$}}
\put(255.00,118.00){\makebox(0,0)[cc]{$_3$}}
\put(180.00,26.00){\makebox(0,0)[cc]{$_{F_1^2}$}}
\put(220.00,26.00){\makebox(0,0)[cc]{$_{C_1^2}$}}
\put(260.00,26.00){\makebox(0,0)[cc]{$_{A_2^2}$}}
\put(180.00,66.00){\makebox(0,0)[cc]{$_{A_1^2}$}}
\put(220.00,66.00){\makebox(0,0)[cc]{$_{F_2^2}$}}
\put(260.00,66.00){\makebox(0,0)[cc]{$_{C_2^2}$}}
\put(180.00,106.00){\makebox(0,0)[cc]{$_{C_0^2}$}}
\put(220.00,106.00){\makebox(0,0)[cc]{$_{A_0^2}$}}
\put(260.00,106.00){\makebox(0,0)[cc]{$_{F_0^2}$}}
\put(165.00,11.00){\makebox(0,0)[cc]{$_{E_2^2}$}}
\put(205.00,11.00){\makebox(0,0)[cc]{$_{E_0^2}$}}
\put(245.00,11.00){\makebox(0,0)[cc]{$_{B_0^2}$}}
\put(165.00,51.00){\makebox(0,0)[cc]{$_{E_1^2}$}}
\put(205.00,51.00){\makebox(0,0)[cc]{$_{D_2^2}$}}
\put(245.00,51.00){\makebox(0,0)[cc]{$_{D_0^2}$}}
\put(165.00,91.00){\makebox(0,0)[cc]{$_{B_2^2}$}}
\put(205.00,91.00){\makebox(0,0)[cc]{$_{D_1^2}$}}
\put(245.00,91.00){\makebox(0,0)[cc]{$_{B_1^2}$}}
\end{picture}
\caption{Toroidal cutouts of $Y_1({\mathcal P})$ and $Y_2({\mathcal P})$}
\end{figure}
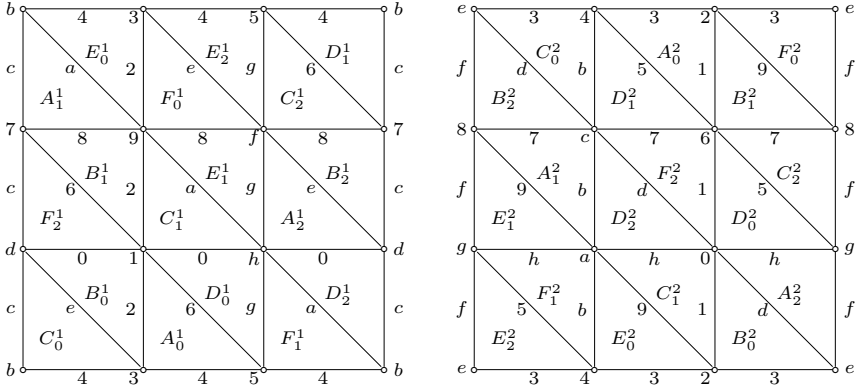

\noindent Both $Y_1({\mathcal P})$ and $Y_2({\mathcal P})$ are
embeddable into a closed orientable surface $T_1$ of genus 1, or
1-torus. Toroidal cutouts of $Y_1({\mathcal P})$ and $Y_2({\mathcal
P})$ are as in Figure 3, which we consider composed by oriented
triangles taken with their orientations derived from those of the
6-cycles of $\mathcal P$ in the proof of Theorem 2, according to the
$2$ alternate operations for $\mathcal P$ mentioned at the end of
Section 3 similar to (2). These oriented copies of $K_3$ are
contractible in $T_1$. They form $2$ collections $\vec{\mathcal
H}_1,\vec{\mathcal H}_2$ of oriented copies $y_i^j$ of $K_3$ closed
under parallel translation, where $y=A,B,C,D,E,F$; $i=0,1,2$ and
$j=1,2$, namely: the 9 oriented triangles of $\vec{\mathcal H}_1$
(resp. $\vec{\mathcal H}_2$) each with horizontal arc below (resp.
above) its opposite vertex. There is also a collection
$\vec{\mathcal H}_0$ of 9 non-contractible oriented triangles in $G$
traceable linearly in 3 different parallel directions, 3 such
triangles per direction, with: {\bf(a)} the orientation of each
participating arc $\vec{e}$ equal to the orientation of the arc of
an oriented triangle in $\vec{\mathcal H}_1$ having the same
end-vertices as $\vec{e}$; {\bf(b)} the arcs of each such oriented
triangle indicated by the (common) middle vertex of the
corresponding 2-arcs in ${\mathcal P}$, as in Section 3. There are
embeddings of $Y_1({\mathcal P})$ and $Y_2({\mathcal P})$ in $T_1$
for which $\vec{\mathcal H}_0$ (resp. $\vec{\mathcal H}_0^{-1}$) and
$\vec{\mathcal H}_1$ (resp. $\vec{\mathcal H}_2$) provide the
composing faces. In addition, each of $\vec{\mathcal H}_1$,
$\vec{\mathcal H}_2$ and $\vec{\mathcal H}_0$ (or $\vec{\mathcal
H}_0^{-1}$) is formed by 3 classes of parallel oriented triangles,
such that any 2 triangles in a class are disjoint. The self-dual
$(9_3)$-configuration in the following theorem is the Pappus $9_3$
\cite{Cox}. Let ${\mathcal H}_i$ be an undirected version of
$\vec{\mathcal H}_i$, for $i=0,1,2$. Let $H_i$ and $\vec{H}_i$ be
respective representatives of ${\mathcal H}_i$ and $\vec{\mathcal
H}_i$, for $i=0,1,2$, and $\vec{H}_0^{-1}$ be a representative of
$\vec{\mathcal H}_0^{-1}$.

\begin{thm}
$Y_1({\mathcal P})$ and $Y_2({\mathcal P})$ are isomorphic tightly
coherent $\{H_0,H_1\cup H_2\}^{P_2}$-homogeneous graphs, as well as
$\{\vec{H}_1,\vec{H}_2\}^{\vec{P}_2}$-,
$\{\vec{H}_1,\vec{H}_0\}^{\vec{P}_2}$- and

\noindent $\{\vec{H}_2,\vec{H}_0^{-1}\}^{\vec{P}_2}$-homogeneous
digraphs. Moreover, each of $Y_1({\mathcal P})$ and $Y_2({\mathcal
P})$ can be taken as the Menger graph of the Pappus self-dual
$(9_3)$-configuration, in $12$ different fashions, by selecting the
point set $\mathcal P$ and the line set ${\mathcal L}\neq{\mathcal
P}$ so that $\{{\mathcal P},{\mathcal L}\}\subset\{V({\mathcal P})$,
${\mathcal H}_0$, ${\mathcal H}_1$, ${\mathcal H}_2\}$ and the
point-line incidence relation either as the inclusion of a vertex in
a triangle or as the containment by a triangle of a vertex or as the
sharing of an edge by $2$ triangles.
\end{thm}

\proof The claimed 12 different forms correspond to the arcs of the
complete graph on vertex set $\{V({\mathcal P})$, ${\mathcal H}_0$,
${\mathcal H}_1$, ${\mathcal H}_2\}$. \qfd\bigskip

\noindent{\bf Acknowledgement.} The author is indebted to the
referee for valuable comments and suggestions.

\end{document}